\documentstyle[12pt,twoside]{article}

\begin{document}

\title{{\bf Poisson-Gradient Dynamical Systems with Bounded Non-Linearity}}
\author{Constantin Udri\c ste, Iulian Duca \\
University Politehnica of Bucharest\\
Department of Mathematics\\
Splaiul Independentei 313\\
060042 Bucharest, ROMANIA\\
Telephone: 4021 4024150\\
e-mail: udriste@mathem.pub.ro}
\date{}
\maketitle

\pagestyle{headings}
\markboth{Poisson-Gradient Dynamical Systems}{Constantin Udri\c ste, Iulian Duca}
\markboth{Poisson-Gradient Dynamical Systems}{Constantin Udri\c ste, Iulian Duca}
\pagestyle{myheadings}

\begin{abstract}
We study the periodical solutions of a Poisson-gradient $PDEs$ system with
bounded non-linearity.

Section 1 introduces the basic spaces and functionals. Section 2 studies the
weak differential of a function and establishes an inequality. Section 3
formulates some conditions under which the action functional is continuously
differentiable. Section 4 analyzes the Poisson-gradient systems and some
conditions that ensure periodical solutions.
\end{abstract}

{\bf Mathematics Subject Classification}: 35J50, 35J55.

{\bf Key words:} variational methods, elliptic systems, periodic solutions,
weak derivatives, potentials.

\section{Introduction}

\noindent We consider the point $T=\left( T^{1},...,T^{p}\right) $ and the parallelepiped $%
T_{0}=\left[ 0,T^{1}\right] \times ...\times \left[ 0,T^{p}\right] $ in
$R^{p}$. We denote by $W_{T}^{1,2}$ the Sobolev space of the functions $%
u\in L^{2}\left[ T_{0},R^{n}\right] $ which have the weak derivative $%
\displaystyle\frac{\partial u}{\partial t}\in L^{2}\left[
T_{0},R^{n}\right] $. The index $T$ from the notation
$W_{T}^{1,2}$ comes from the fact that the weak derivatives are
defined using the space $C_{T}^{\infty }$ of all indefinitely
differentiable multiple $T$-periodic functions from $R^{p}$ into
$R^{n}$. We denote by $H_{T}^{1}$ the Hilbert space $W_{T}^{1,2}.$
The norm used in $H_{T}^{1}$ is the one induced by the scalar
product
\[
\left\langle u,v\right\rangle =\int_{T_{0}}\left( \delta _{ij}u^{i}\left(
t\right) v^{j}\left( t\right) +\delta _{ij}\delta ^{\alpha \beta }\frac{%
\partial u^{i}}{\partial t^{\alpha }}\left( t\right) \frac{\partial v^{j}}{%
\partial t^{\beta }}\left( t\right) \right) dt^{1}\wedge ...\wedge dt^{p}.
\]
These are induced by the scalar product (Riemaniann metric)
\[
G=\left(
\begin{array}{cc}
\delta _{ij} & 0 \\
0 & \delta ^{\alpha \beta }\delta _{ij}
\end{array}
\right)
\]
on $R^{n+np}$ (multiphase space) and\ its
associated Euclidean norm. We shall also use the scalar product $\left( u,v\right) =\delta
_{ij}u^{i}v^{j}$ and the norm $\left| u\right| =\sqrt{\delta _{ij}u^{i}u^{j}}
$ simultaneously, from the Euclidean space $R^{n}$.

Let $t=\left(
t^{1},...,t^{p}\right) $ be a generic point in $R^{p}$. Then the opposite
faces of the parallelepiped $T_{0}$ can be described by the equations
\[
S_{i}^{-}:t^{i}=0,S_{i}^{+}:t^{i}=T^{i}
\]
for each $i=1,...,p.$ We shall study the minimum of the action
\[
\varphi \left( u\right) =\displaystyle\int_{T_{0}}L\left( t,u\left( t\right)
,\frac{\partial u}{\partial t}\right) dt^{1}\wedge ...\wedge dt^{p},
\]
\[
L\left( t,u\left( t\right) ,\displaystyle\frac{\partial u}{\partial t}%
\right) =\displaystyle\frac{1}{2}\left| \displaystyle\frac{\partial u}{%
\partial t}\right| ^{2}+F\left( t,u\left( t\right) \right)
\]
on the space $H_{T}^{1}$, considering that the potential function
$F$ has the property of bounded non-linearity. We use the method
of the minimizing sequences and the coercitivity condition
$\displaystyle\int_{T_{0}}F\left( t,u\left( t\right) \right)
dt^{1}\wedge ...\wedge dt^{p}\;\rightarrow \;\infty $ when $\left|
u\right| \;\rightarrow \;\infty .$ The extremals of the action
$\varphi $ verifies the Euler-Lagrange equations with the boundary
conditions
\[
u\mid _{S_{i}^{-}}=u\mid _{S_{i}^{+}},\frac{\partial u}{\partial t}\mid
_{S_{i}^{-}}=\frac{\partial u}{\partial t}\mid _{S_{i}^{+}},i=1,...,p.
\]
Due to the particularity of the Lagrangian $L$, the Euler-Lagrange equations
reduce to a $PDEs$ system of the Poisson-gradient type
\[
\Delta u\left( t\right) =\nabla F\left( t,u\left( t\right) \right) .
\]

The aim of this paper is to discuss the existence of solutions of this PDEs system with suitable boundary conditions.
More precisely, we extend the theory in [2] from single-time to multi-time field theory, developing the ideas in the papers [6], [7], [9].
 In this way we find positive answers for the existence of multi-periodical solutions of Euler-Lagrange equations that are Poisson-gradient PDEs with bounded non-linearity. The results can be applied to the multi-time geometric dynamics ([5], [8], [10]-[12]).

\section{On the weak differential of a function}

We consider $C_{T}^{\infty }$ the space of the indefinitely differentiable
functions multiple periodical with the period $T=\left(
T^{1},...,T^{p}\right) $, defined on $R^{p}$ taking values in $R^{n}$. We
know that $C_{T}^{\infty }\subset W_{T}^{1,2}$. We establish some conditions
satisfied by a function $u\in L^{1}\left[ T_{0},R^{n}\right] $ which has a
weak differential.

{\bf Theorem 1.} {\it Let $u,v_{\alpha }\in L^{1}\left[ T_{0},R^{n}\right]
,\alpha =1,...,p$, such that $v_{\alpha }dt^{\alpha }=\left( v_{\alpha
}^{1}dt^{\alpha },...,v_{\alpha }^{1}dt^{\alpha }\right) $ is an integrable
vector form. We consider $\stackrel{\frown }{OT}$ an arbitrary curve\ from $%
T_{0},$ having the endings at $O=\left( 0,...,0\right) $ and $T=\left(
T^{1},...,T^{p}\right) .$ }

{\it If
$$
\int_{\stackrel{\frown }{OT}}\left( u,df\right) =-\int_{\stackrel{\frown }{OT%
}}\left( v_{\alpha }dt^{\alpha },f\right) ,\eqno(1)
$$
for any $f\in C_{T}^{\infty }$, then $\displaystyle\int_{\stackrel{\frown }{%
OT}}v_{\alpha }dt^{\alpha }=0$ and it exists $c\in R^{n}$ such that $u\left(
t\right) =\displaystyle\int_{\stackrel{\frown }{Ot}}v_{\alpha }ds^{\alpha
}+c $. Also $u\left( 0\right) =u\left( T\right) .$}

{\bf Proof}. We choose $f=e^{i}=\left( 0,...,0,1,0,...,0\right) $, with the
value 1 on the position $i$. From the relation $\left( 1\right) $ we have $%
0=-\displaystyle\int_{\stackrel{\frown }{OT}}v_{\alpha }^{i}dt^{\alpha }$
and hence $\displaystyle\int_{\stackrel{\frown }{OT}}v_{\alpha }dt^{\alpha
}=0$.

We define $w\in C\left( T_{0},R^{n}\right) $ by $w\left( t\right) =$ $%
\displaystyle\int_{\stackrel{\frown }{Ot}}v_{\alpha }ds^{\alpha },t\in
\stackrel{\frown }{OT}$. By Fubini Theorem, the function $w$ satisfies the
relation
\[
\displaystyle\int_{\stackrel{\frown }{OT}}\left( w\left( t\right) ,df\right)
=\displaystyle\int_{\stackrel{\frown }{OT}}\left( \displaystyle\int_{%
\stackrel{\frown }{Ot}}v_{\alpha }ds^{\alpha },df\right) =\displaystyle\int_{%
\stackrel{\frown }{OT}}(\displaystyle\int_{\stackrel{\frown }{sT}}\left(
v_{\alpha },df\right) ds^{\alpha })
\]
\[
=\displaystyle\int_{\stackrel{\frown }{OT}}\left( v_{\alpha },f\left(
T\right) -f\left( s\right) \right) ds^{\alpha }=-\int_{\stackrel{\frown }{OT}%
}\left( v_{\alpha },f\left( s\right) \right) ds^{\alpha }=\int_{\stackrel{%
\frown }{OT}}\left( u,df\right) .
\]
This means that
$$
\int_{\stackrel{\frown }{OT}}\left( u-w,df\right) =0.\eqno(2)
$$
We consider now $\gamma :\left[ a,b\right] \rightarrow
T_{0},\gamma \left( \xi \right) =\left( t^{1}\left( \xi \right)
,...,t^{p}\left( \xi \right) \right) ,\gamma \left( a\right) =O,$
$\ \gamma \left( b\right) =T$, a parameterization of the curve
$\stackrel{\frown }{OT}$. The equality $\left( 2\right) $ becomes
\[
\int_{a}^{b}\left( u\left( t\left( \xi \right) \right) -w\left( t\left( \xi
\right) \right) ,\left( \frac{\partial f^{1}}{\partial t^{\alpha }}\frac{%
dt^{\alpha }}{d\xi },...,\frac{\partial f^{n}}{\partial t^{\alpha }}\frac{%
dt^{\alpha }}{d\xi }\right) \right) d\xi =0,
\]
for any $f\in C_{T}^{\infty }$. We will particularize for the function
sequences
\[
f_{j}^{\left( k\right) }\left( t\right) =\left\{
\begin{array}{c}
\cos \\
\sin
\end{array}
\right\} \left( \frac{2k\pi t^{j}}{T^{j}}\right) e^{j},\;k\in N\setminus
\{0\},\;1\leq j\leq n
\]
and we observe that (see the Fourier series theory) $u\left( t\right)
-w\left( t\right) =c$, $c\in R^{n}$ almost everywhere in $T_{0}$ $($the
constant is the only function orthogonal to the previous sequences). By
replacing $w\left( t\right) $, we find that $u\left( t\right) =\displaystyle%
\int_{\stackrel{\frown }{Ot}}v_{\alpha }ds^{\alpha }+c$ for any $t\in
\stackrel{\frown }{OT}.$ The function $u$ satisfies $u\left( 0\right) =c$
and $u\left( T\right) =\displaystyle\int_{\stackrel{\frown }{OT}}v_{\alpha
}ds^{\alpha }+c=c,$ so $u\left( 0\right) =u\left( T\right) .$ On the other
side, the relation $u\left( t\right) -u\left( \tau \right) =\displaystyle%
\int_{\stackrel{\frown }{\tau t}}v_{\alpha }ds^{\alpha }$ implies that $%
u\left( t\right) =\displaystyle\int_{\stackrel{\frown }{\tau t}}v_{\alpha
}ds^{\alpha }+u\left( \tau \right) $. The 1-form $v_{\alpha }dt^{\alpha }$
is called {\it weak differential }of the function $u$. By a Fourier series
argument, the weak differential, if it exists, is unique. The weak
differential of $u$ will be denoted by $du$. The existence of $du$ implies $%
u\left( 0\right) =u\left( T\right) .$

{\bf Theorem 2.} {\it If $u=\left( u^{1},...,u^{n}\right) \in L^{2}\left[
T_{0},R^{n}\right] ,$ $\ \left| u\left( t\right) \right|
^{2}=\delta _{ij}u^{i}\left( t\right) u^{j}\left( t\right) $, then}
\[
\left| \int_{T_{0}}u\left( t\right) dt^{1}\wedge ...\wedge dt^{p}\right|
\leq \left( nT^{1}...T^{p}\right) ^{\frac{1}{2}}\left( \int_{T_{0}}\left|
u\left( t\right) \right| ^{2}dt^{1}\wedge ...\wedge dt^{p}\right) ^{\frac{1}{%
2}}.
\]

{\bf Proof}. Successively we have the relations
\[
\left| \int_{T_{0}}u\left( t\right) dt^{1}\wedge ...\wedge dt^{p}\right|
=\left| \int_{T_{0}}\left( u^{1}\left( t\right) ,...,u^{n}\left( t\right)
\right) dt^{1}\wedge ...\wedge dt^{p}\right|
\]
\[
=\left| \left( \int_{T_{0}}u^{1}\left( t\right) dt^{1}\wedge ...\wedge
dt^{p},...,\int_{T_{0}}u^{n}\left( t\right) dt^{1}\wedge ...\wedge
dt^{p}\right) \right|
\]
\[
=\left( \left( \int_{T_{0}}u^{1}\left( t\right) dt^{1}\wedge ...\wedge
dt^{p}\right) ^{2}+...+\left( \int_{T_{0}}u^{1}\left( t\right) dt^{1}\wedge
...\wedge dt^{p}\right) ^{2}\right) ^{\frac{1}{2}}
\]
\[
\leq \left| \int_{T_{0}}u^{1}\left( t\right) dt^{1}\wedge ...\wedge
dt^{p}\right| +...+\left| \int_{T_{0}}u^{n}\left( t\right) dt^{1}\wedge
...\wedge dt^{p}\right|
\]
\[
\leq \int_{T_{0}}\left( \left| u^{1}\left( t\right) \right| +...+\left|
u^{n}\left( t\right) \right| \right) dt^{1}\wedge ...\wedge dt^{p}
\]
\[
=\int_{T_{0}}\left( \left( \left| u^{1}\left( t\right) \right| ,...,\left|
u^{n}\left( t\right) \right| \right) ,\left( 1,...,1\right) \right)
dt^{1}\wedge ...\wedge dt^{p}.\]
Using the Cauchy-Schwartz inequality, we obtain
\[
\left| \int_{T_{0}}u\left( t\right) dt^{1}\wedge ...\wedge dt^{p}\right|
\]
\[
\leq \left( \int_{T_{0}}\left( \left| u^{1}\left( t\right) \right|
^{2}+...+\left| u^{n}\left( t\right) \right| ^{2}\right) dt^{1}\wedge
...\wedge dt^{p}\right) ^{\frac{1}{2}}\left( \int_{T_{0}}ndt^{1}\wedge
...\wedge dt^{p}\right) ^{\frac{1}{2}}
\]
\[
=\left( nT^{1}...T^{p}\right) ^{\frac{1}{2}}\left( \int_{T_{0}}\left|
u\left( t\right) \right| ^{2}dt^{1}\wedge ...\wedge dt^{p}\right) ^{\frac{1}{%
2}}.
\]

\section{Continuously differentiable action}

The next theorem establishes some conditions in which the action
\[
\varphi :W_{T}^{1,2}\rightarrow R,\varphi \left( u\right)
=\int_{T_{0}}L\left( t,u\left( t\right) ,\frac{\partial u}{\partial t}\left(
t\right) \right) dt^{1}\wedge ...\wedge dt^{p}
\]
is continuously differentiable. In this way we extend the particular case $p=1$, studied in
[3, Theorem 1.4].

{\bf Theorem 3}. {\it We consider $L:T_{0}\times R^{n}\times
R^{np}\rightarrow R,\left( t,x,y\right) \rightarrow L\left( t,x,y\right) $,
a measurable function in $t$ for any $\left( x,y\right) \in R^{n}\times
R^{np}$ and with the continuous partial derivatives in x and y for any $t\in
T_{0}$. If here exist $a\in C^{1}\left( R^{+},R^{+}\right) $ with the
derivative $a^{\prime }$\ bounded from above, $b\in C\left(
T_{0},R^{n}\right) $ such that for any $t\in T_{0}$ and any $\left(
x,y\right) \in R^{n}\times R^{np}$ to have }
$$
\begin{array}{l}
\left| L\left( t,x,y\right) \right| \leq a\left( \left| x\right| +\left|
y\right| ^{2}\right) b\left( t\right) , \\
\noalign{\medskip}\left| \nabla _{x}L\left( t,x,y\right) \right| \leq
a\left( \left| x\right| \right) b\left( t\right) , \\
\noalign{\medskip}\left| \nabla _{y}L\left( t,x,y\right) \right| \leq
a\left( \left| y\right| \right) b\left( t\right) ,
\end{array}
\eqno(3)
$$
{\it then, the functional $\varphi $ has continuous partial derivatives in $%
W_{T}^{1,2}$ and his gradient derives from the formula }
$$
\begin{array}{lcl}
\left( \nabla \varphi \left( u\right) ,v\right) & = & \displaystyle%
\int_{T_{0}}\left[ \left( \nabla _{x}L\left( t,u\left( t\right) ,%
\displaystyle\frac{\partial u}{\partial t}\right) ,v\left( t\right) \right)
\right. \\
\noalign{\medskip} & + & \left. \left( \nabla _{y}L\left( t,u\left( t\right)
,\displaystyle\frac{\partial u}{\partial t}\left( t\right) \right) ,%
\displaystyle\frac{\partial v}{\partial t}\left( t\right) \right) \right]
dt^{1}\wedge ...\wedge dt^{p}.
\end{array}
\eqno(4)
$$

{\bf Proof}. It is enough to prove that $\varphi $ has the derivative $%
\varphi ^{\prime }\left( u\right) \in \left( W_{T}^{1,2}\right) ^{\ast }$
given by the relation $\left( 4\right) $ and the function $\varphi ^{\prime
}:W_{T}^{1,2}\rightarrow \left( W_{T}^{1,2}\right) ^{\ast }$, $u\rightarrow
\varphi ^{\prime }\left( u\right) $ is continuous. We consider $u,v\in
W_{T}^{1,2}$, $t\in T_{0}$, $\lambda \in \left[ -1,1\right] $. We build the
functions
\[
F\left( \lambda ,t\right) =L\left( t,u\left( t\right) +\lambda v\left(
t\right) ,\frac{\partial u}{\partial t}\left( t\right) +\lambda \frac{%
\partial v}{\partial t}\left( t\right) \right)
\]
and
\[
\Psi \left( \lambda \right) =\int_{T_{0}}F\left( \lambda ,t\right)
dt^{1}\wedge ...\wedge dt^{p}.
\]
Because the derivative $a^{\prime }$ is bounded from above, exist\ $M>0$
such that $\displaystyle\frac{a\left( \left| u\right| \right) -a\left(
0\right) }{\left| u\right| }=a^{\prime }\left( c\right) \leq M.$ This means
that $a\left( \left| u\right| \right) \leq M\left| u\right| +a\left(
0\right) .$ On the other side

\[
\frac{\partial F}{\partial \lambda }\left( \lambda ,t\right) =\left( \nabla
_{x}L\left( t,u\left( t\right) +\lambda v\left( t\right) ,\frac{\partial u}{%
\partial t}\left( t\right) +\lambda \frac{\partial v}{\partial t}\left(
t\right) \right) ,v\left( t\right) \right)
\]
\[
+\left( \nabla _{y}L\left( t,u\left( t\right) +\lambda v\left( t\right) ,%
\frac{\partial u}{\partial t}\left( t\right) +\lambda \frac{\partial v}{%
\partial t}\left( t\right) \right) ,\frac{\partial v}{\partial t}\left(
t\right) \right) \leq a\left( \left| u\left( t\right) +\lambda v\left(
t\right) \right| \right)
\]
\[
b\left( t\right) \left| v\left( t\right) \right| +a\left( \left| \frac{%
\partial u}{\partial t}\left( t\right) +\lambda \frac{\partial v}{\partial t}%
\left( t\right) \right| \right) b\left( t\right) \left| \frac{\partial v}{%
\partial t}\left( t\right) \right|
\]
\[
\leq b_{0}\left( M\left( \left| u\left( t\right) \right| +\left| v\left(
t\right) \right| \right) +a\left( 0\right) \right) \left| v\left( t\right)
\right| +
\]
\[
b_{0}\left( M\left( \left| \frac{\partial u}{\partial t}\left(
t\right) \right| +\left| \frac{\partial v}{\partial t}\left( t\right)
\right| \right) +a\left( 0\right) \right) \left| \frac{\partial v}{\partial t%
}\left( t\right) \right| ,
\]
where
\[
b_{0}=\displaystyle\max_{t\in T_{0}}b\left( t\right) .
\]
Then, we have $\left| \displaystyle\frac{\partial F}{\partial \lambda }%
\left( \lambda ,t\right) \right| \leq d\left( t\right) \in L^{1}\left(
T_{0},R^{+}\right) $. Then Leibniz formula of differentiation under integral
sign is applicable and
\[
\frac{\partial \Psi }{\partial \lambda }\left( 0\right) =\int_{T_{0}}\frac{%
\partial F}{\partial \lambda }\left( 0,t\right) dt^{1}\wedge ...\wedge
dt^{p}=\int_{T_{0}}\left[ \left( \nabla _{x}L\left( t,u\left( t\right) ,%
\frac{\partial u}{\partial t}\left( t\right) \right) ,v\left( t\right)
\right) \right.
\]
\[
\left. +\left( \nabla _{y}L\left( t,u\left( t\right) ,\frac{\partial u}{%
\partial t}\left( t\right) \right) ,\frac{\partial v}{\partial t}\left(
t\right) \right) \right] dt^{1}\wedge ...\wedge dt^{p}.
\]
Moreover,
\[
\left| \nabla _{x}L\left( t,u\left( t\right) ,\frac{\partial u}{\partial t}%
\left( t\right) \right) \right| \leq b_{0}\left( M\left| u\left( t\right)
\right| +\left| a\left( 0\right) \right| \right) \in L^{1}\left(
T_{0},R^{+}\right)
\]
and
\[
\left| \nabla _{y}L\left( t,u\left( t\right) ,\frac{\partial u}{\partial t}%
\left( t\right) \right) \right| \leq b_{0}\left( M\left| \frac{\partial u}{%
\partial t}\left( t\right) \right| +\left| a\left( 0\right) \right| \right)
\in L^{2}\left( T_{0},R^{+}\right) .
\]
That is why

\[
\int_{T_{0}}\left[ \left( \nabla _{x}L\left( t,u\left( t\right) ,\frac{%
\partial u}{\partial t}\left( t\right) \right) ,v\left( t\right) \right)
\right.
\]
\[
\left. +\left( \nabla _{y}L\left( t,u\left( t\right) ,\frac{\partial u}{%
\partial t}\left( t\right) \right) ,\frac{\partial v}{\partial t}\left(
t\right) \right) \right] dt^{1}\wedge ...\wedge dt^{p}
\]

\[
\leq \int_{T_{0}}\left| \nabla _{x}L\left( t,u\left( t\right) ,\frac{%
\partial u}{\partial t}\left( t\right) \right) \right| \left| v\left(
t\right) \right| dt^{1}\wedge ...\wedge dt^{p}
\]
\[
+\int_{T_{0}}\left| \nabla _{y}L\left( t,u\left( t\right) ,\frac{\partial u}{%
\partial t}\left( t\right) \right) \right| \left| \frac{\partial v}{\partial
t}\left( t\right) \right| dt^{1}\wedge ...\wedge dt^{p}
\]
\[
\leq b_{0}\int_{T_{0}}\left( M\left| u\left( t\right) \right| +\left|
a\left( 0\right) \right| \right) \left| v\left( t\right) \right|
dt^{1}\wedge ...\wedge dt^{p}
\]
\[
+b_{0}\int_{T_{0}}\left( M\left| \frac{\partial u}{\partial t}\left(
t\right) \right| +\left| a\left( 0\right) \right| \right) \left| \frac{%
\partial v}{\partial t}\left( t\right) \right| dt^{1}\wedge ...\wedge dt^{p}.\]
By using the inequality Cauchy-Schwartz, we find
\[
\begin{array}{l}
\left| \displaystyle\frac{\partial \Psi }{\partial \lambda }\left( 0\right)
\right| \leq b_{0}\left( \displaystyle\int_{T_{0}}\left( M\left| u\left(
t\right) \right| +\left| a\left( 0\right) \right| \right) ^{2}dt^{1}\wedge
...\wedge dt^{p}\right) ^{\frac{1}{2}} \\
\noalign{\medskip}\left( \displaystyle\int_{T_{0}}\left| v\left(
t\right) \right| ^{2}dt^{1}\wedge ...\wedge dt^{p}\right) ^{\frac{1}{2}} \\
\noalign{\medskip} + b_{0}\left( \displaystyle\int_{T_{0}}\left( M\left| %
\displaystyle\frac{\partial u}{\partial t}\left( t\right) \right| +\left|
a\left( 0\right) \right| \right) ^{2}dt^{1}\wedge ...\wedge dt^{p}\right) ^{%
\frac{1}{2}}\\

\left( \displaystyle\int_{T_{0}}\left| \displaystyle\frac{%
\partial v}{\partial t}\left( t\right) \right| ^{2}dt^{1}\wedge ...\wedge
dt^{p}\right) ^{\frac{1}{2}}
\end{array}
\]

\[
\begin{array}{l}
\noalign{\medskip}\leq C_{1}\left( \displaystyle\int_{T_{0}}\left| v\left(
t\right) \right| ^{2}dt^{1}\wedge ...\wedge dt^{p}\right) ^{\frac{1}{2}%
}+C_{2}\left( \displaystyle\int_{T_{0}}\left| \displaystyle\frac{\partial v}{%
\partial t}\left( t\right) \right| ^{2}dt^{1}\wedge ...\wedge dt^{p}\right)
^{\frac{1}{2}} \\
\noalign{\medskip}\leq \max \left\{ C_{1},C_{2}\right\} 2^{\frac{1}{2}%
}\left( \displaystyle\int_{T_{0}}\left( \left| v\left( t\right) \right|
^{2}+\left| \displaystyle\frac{\partial v}{\partial t}\left( t\right)
\right| ^{2}\right) dt^{1}\wedge ...\wedge dt^{p}\right) ^{\frac{1}{2}%
}=C\left\| v\right\| .
\end{array}
\]
By consequence, the action $\varphi $ has the derivative $\varphi ^{\prime
}\in \left( W_{T}^{1,2}\right) ^{\ast }$ given by $\left( 4\right) $. The
Krasnoselski theorem and the hypothesis $\left( 3\right) $ imply the fact
that the application $u\rightarrow \left( \nabla _{x}L\left( \cdot ,u,%
\displaystyle\frac{\partial u}{\partial t}\right) ,\nabla _{y}L\left( \cdot
,u,\displaystyle\frac{\partial u}{\partial t}\right) \right) ,$ from $%
W_{T}^{1,2}$ to $L^{1}\times L^{2},$ is continuous, so $\varphi ^{\prime }$
is continuous from $W_{T}^{1,2}$ to $\left( W_{T}^{1,2}\right) ^{\ast }$ and
the proof is complete.

\section{Poisson-gradient systems and their periodical solutions}

\subsection{Multi-time Euler-Lagrange equations}

We consider the multi-time variable $t=\left( t^{1},...,t^{p}\right) \in
R^{p}$, the functions $x^{i}:R^{p}\rightarrow R,\left( t^{1},...,t^{p}\right)
\rightarrow $ $x^{i}\left( t^{1},...,t^{p}\right) ,i=1,...n,$ and the
partial velocities $x_{\alpha }^{i}=\displaystyle\frac{\partial x^{i}}{%
\partial t^{\alpha }},\alpha =1,...,p.$ The Lagrangian
\[
L:R^{p+n+np}\rightarrow R,\left( t^{\alpha },x^{i},x_{\alpha }^{i}\right)
\rightarrow L\left( t^{\alpha },x^{i},x_{\alpha }^{i}\right)
\]
determines the Euler-Lagrange equations
\[
\frac{\partial }{\partial t^{\alpha }}\frac{\partial L}{\partial x_{\alpha
}^{i}}=\frac{\partial L}{\partial x^{i}},\;i=1,...,n,\;\alpha =1,...p
\]
($PDEs$ system of second order in the n-dimensional space). We remark that
in the left hand member we have summation after the index $\alpha $ (trace).

\subsection{An action that produces Poisson-gradient systems}

Let $\alpha =1,...,p,i=1,...,n,u^{i}:T_{0}\rightarrow R,t=\left(
t^{1},...,t^{p}\right) \rightarrow u^{i}\left( t^{1},...,t^{p}\right)
,u:T_{0}\rightarrow R^{n},u\left( t\right) =\left( u^{1}\left( t\right)
,...,u^{n}\left( t\right) \right) ,\displaystyle u_{\alpha }^{i}=\frac{%
\partial u^{i}}{\partial t^{\alpha }},\displaystyle\frac{\partial u}{%
\partial t}=\left( u_{\alpha }^{i}\right) .$

We consider the Lagrangian
\[
L:T_{0}\times R^{n}\times R^{np}\rightarrow R,\left( t^{\alpha
},u^{i},u_{\alpha }^{i}\right) \rightarrow L\left( t^{\alpha
},u^{i},u_{\alpha }^{i}\right) ,
\]
\[
L\left( t^{\alpha },u^{i},u_{\alpha }^{i}\right) =\frac{1}{2}\left| \frac{%
\partial u}{\partial t}\right| ^{2}+F(t,u\left( t\right) ).
\]
A function $u$ (field) that realizes the minimum of the action
\[
\varphi \left( u\right) =\displaystyle\int_{T_{0}}L\left( t,u\left( t\right)
,\displaystyle\frac{\partial u}{\partial t}\left( t\right) \right)
dt^{1}\wedge ...\wedge dt^{p},
\]
verifies a $PDEs$ system of Poisson-gradient type (Euler-Lagrange equations
on $H_{T}^{1}$)\
\[
\Delta u\left( t\right) =\nabla F\left( t,u\left( t\right) \right) ,
\]
together with the boundary conditions
\[
u\mid _{S_{i}^{-}}=u\mid _{S_{i}^{+}},\frac{\partial u}{\partial t}\mid
_{S_{i}^{-}}=\frac{\partial u}{\partial t}\mid _{S_{i}^{+}},i=1,...,p.
\]

\subsection{Periodical solutions of Poisson-gradient dynamical systems with
bounded non-linearity}

{\bf Theorem 4}. {\it Suppose the function$\ F:T_{0}\times R^{n}\rightarrow
R,\left( t,u\right) \rightarrow F\left( t,u\right) $ satisfies four
properties: }

{\it 1) $F\left( t,u\right) $ is measurable in $t$ for any $u\in R^{n}$ and
it is continuously differentiable in $u$ for any $t\in T_{0}$}${\it ,}${\it %
\ }

{\it 2) There exist the functions $a\in C^{1}\left( R^{+},R^{+}\right) $
with the derivative $a^{\prime }$\ bounded from above and $b\in C\left(
T_{0},R^{+}\right) $ such that for any $t\in T_{0}$ and any $u\in R^{n}$ to
have $\left| F\left( t,u\right) \right| \leq a\left( \left| u\right| \right)
b\left( t\right) $ and $\left| \nabla _{u}F\left( t,u\right) \right| \leq
a\left( \left| u\right| \right) b\left( t\right) $}${\it ,}${\it \ }

{\it 3) It exists $g\in C^{1}\left( T_{0},R\right) $ such that for any $t\in
T_{0}$ and any $u\in R^{n}$, to have
\[
\left| \nabla _{u}F\left( t,u\right) \right| \leq g\left( t\right) .
\]
}

\bigskip 4){\it \ The action $\varphi _{1}\left( u\right) =\displaystyle%
\int_{T_{0}}F\left( t,u\left( t\right) \right) dt^{1}\wedge
...\wedge dt^{p}$ is weakly lower semi-continuous.}

{\it If $\displaystyle\int_{T_{0}}F\left( t,u\right) dt^{1}\wedge ...\wedge
dt^{p}\rightarrow \infty $ when $\left| u\right| \rightarrow \infty $}${\it ,%
}${\it \ then the Dirichlet problem }
\[
{\it \Delta u\left( t\right) =\nabla F\left( t,u\left( t\right) \right) ,}
\]
{\it \ \ }
\[
{\it u\mid _{S_{i}^{-}}=u\mid _{S_{i}^{+}},\frac{\partial u}{\partial t}\mid
_{S_{i}^{-}}=\frac{\partial u}{\partial t}\mid _{S_{i}^{+}},i=1,...,p,}
\]
{\it \ has at least a solution which minimizes the action}
\[
{\it \varphi \left( u\right) =\displaystyle\int_{T_{0}}\left[ \displaystyle%
\frac{1}{2}\left| \displaystyle\frac{\partial u}{\partial t}\right|
^{2}+F\left( t,u\left( t\right) \right) \right] dt^{1}\wedge ...\wedge dt^{p}%
}
\]
{\it in} {\it $H_{T}^{1}$}.

{\bf Proof}. We consider $u=\overline{u}+\widetilde{u}$, where $\overline{u}=%
 \displaystyle\frac{1}{T^1...T^p}\displaystyle\int_{T_{0}}u\left( t\right) dt^{1}\wedge ...\wedge dt^{p}$.
Then
\[
\varphi \left( u\right) =\int_{T_{0}}\left[ \frac{1}{2}\left| \frac{\partial
u}{\partial t}\right| ^{2}+F\left( t,u\left( t\right) \right) \right]
dt^{1}\wedge ...\wedge dt^{p}
\]
\[
=\int_{T_{0}}\left[ \frac{1}{2}\left| \frac{\partial u}{\partial t}\right|
^{2}+F\left( t,\overline{u}\right) -F\left( t,\overline{u}\right) +F\left(
t,u\left( t\right) \right) \right] dt^{1}\wedge ...\wedge dt^{p}
\]
\[
=\int_{T_{0}}\left[ \frac{1}{2}\left| \frac{\partial u}{\partial t}\right|
^{2}+F\left( t,\overline{u}\left( t\right) \right) \right] dt^{1}\wedge
...\wedge dt^{p}
\]
\[
+\int_{T_{0}}\int_{0}^{1}\left( \nabla _{u}F\left( t,\overline{u}+s%
\widetilde{u}\left( t\right) \right) ,\widetilde{u}\left( t\right) \right)
ds\wedge dt^{1}\wedge ...\wedge dt^{p}.
\]
According to property 3) from the hypothesis, we have the inequality
\[
\left( \nabla _{u}F\left( t,\overline{u}+s\widetilde{u}\left( t\right)
\right) ,\widetilde{u}\left( t\right) \right) \leq \left| \nabla _{u}F\left(
t,\overline{u}+s\widetilde{u}\left( t\right) \right) \right| \left|
\widetilde{u}\left( t\right) \right| \leq \left| g\left( t\right) \right|
\left| \widetilde{u}\left( t\right) \right|
\]
from which we obtain the relation
\[
-\left| g\left( t\right) \right| \left| \widetilde{u}\left( t\right) \right|
\leq \left( \nabla _{u}F\left( t,\overline{u}+s\widetilde{u}\left( t\right)
\right) ,\widetilde{u}\left( t\right) \right)
\]
for any $t\in T_{0}$. By using this inequality we obtain
\[
\varphi \left( u\right) =\int_{T_{0}}\frac{1}{2}\left| \frac{\partial u}{%
\partial t}\right| ^{2}dt^{1}\wedge ...\wedge dt^{p}+\int_{T_{0}}F\left( t,%
\overline{u}\right) dt^{1}\wedge ...\wedge dt^{p}
\]
\[
-\int_{T_{0}}\left| g\left( t\right) \right| \left| \widetilde{u}\left(
t\right) \right| dt^{1}\wedge ...\wedge dt^{p}
\]
\[
\geq \int_{T_{0}}\frac{1}{2}\left| \frac{\partial u}{\partial t}\right|
^{2}dt^{1}\wedge ...\wedge dt^{p}+\int_{T_{0}}F\left( t,\overline{u}\right)
dt^{1}\wedge ...\wedge dt^{p}
\]
\[
-g_{0}\int_{T_{0}}\left| \widetilde{u}\left( t\right) \right| dt^{1}\wedge
...\wedge dt^{p},
\]
where $g_{0}=\displaystyle\max_{t\in T_{0}}\left| g\left( t\right) \right| $.
According to the multi-time Wirtinger inequality [9], it exists $C_1>0$ such
that
\[
\int_{T_{0}}\left| \widetilde{u}\left( t\right) \right| dt^{1}\wedge
...\wedge dt^{p}\leq C_{1}\left( \int_{T_{0}}\left| \frac{\partial u}{%
\partial t}\left( t\right) \right| ^{2}dt^{1}\wedge ...\wedge dt^{p}\right)
^{\frac{1}{2}}.
\]
This means that
\[
\varphi \left( u\right) \geq \int_{T_{0}}\frac{1}{2}\left| \frac{\partial u}{%
\partial t}\right| ^{2}dt^{1}\wedge ...\wedge dt^{p}+\int_{T_{0}}F\left( t,%
\overline{u}\right) dt^{1}\wedge ...\wedge dt^{p}
\]
\[
-g_{0}C_{1}\left( \int_{T_{0}}\left| \frac{\partial u}{\partial t}\left(
t\right) \right| ^{2}dt^{1}\wedge ...\wedge dt^{p}\right) ^{\frac{1}{2}}.
\]
Of course, if $\left\| u\right\| \rightarrow \infty $, then, from the
relation $\left\| u\right\| \leq \left\| \overline{u}\right\| +\left\|
\widetilde{u}\right\| $ it follows that $\left\| \overline{u}\right\|
\rightarrow \infty $ or $\left\| \widetilde{u}\right\| \rightarrow \infty $.
Because $\overline{u}$ is constant in $R^{n}$, we have the equalities
\[
\left\| \overline{u}\right\| =\left\| \overline{u}\right\|
_{W_{T}^{1,2}}=\left( \int_{T_{0}}\left( \left| \overline{u}\right|
^{2}+\left| \frac{\partial \overline{u}}{\partial t}\right| \right)
^{2}dt^{1}\wedge ...\wedge dt^{p}\right) ^{\frac{1}{2}}
\]
\[
=\left( \int_{T_{0}}\left| \overline{u}\right| ^{2}dt^{1}\wedge ...\wedge
dt^{p}\right) ^{\frac{1}{2}}=\left| \overline{u}\right| \left(
T^{1}...T^{p}\right) ^{\frac{1}{2}}.
\]
This means that if $\left\| \overline{u}\right\| \rightarrow \infty ,$ then $%
\left| \overline{u}\right| \rightarrow \infty $. Consequently using the
hypothesis, we obtain
$$
\int_{T_{0}}F\left( t,\overline{u}\right) dt^{1}\wedge ...\wedge
dt^{p}\rightarrow \infty .\eqno(5)
$$
Also
\[
\left\| \widetilde{u}\right\| =\left( \int_{T_{0}}\left( \left| \widetilde{u}%
\left( t\right) \right| ^{2}+\left| \frac{\partial \widetilde{u}}{\partial t}%
\left( t\right) \right| \right) ^{2}dt^{1}\wedge ...\wedge dt^{p}\right) ^{%
\frac{1}{2}}
\]
\[
=\left( \int_{T_{0}}\left( \left| \widetilde{u}\left( t\right) \right|
^{2}+\left| \frac{\partial u}{\partial t}\left( t\right) \right| \right)
^{2}dt^{1}\wedge ...\wedge dt^{p}\right) ^{\frac{1}{2}}.
\]
With the Wirtinger inequality we obtain
\[
\left\| \widetilde{u}\right\| \leq \left( \int_{T_{0}}\left( C\left| \frac{%
\partial \widetilde{u}}{\partial t}\left( t\right) \right| ^{2}+\left| \frac{%
\partial u}{\partial t}\left( t\right) \right| \right) ^{2}dt^{1}\wedge
...\wedge dt^{p}\right) ^{\frac{1}{2}}
\]
\[
=\left( C+1\right) \left( \int_{T_{0}}\left| \frac{\partial u}{\partial t}%
\left( t\right) \right| ^{2}dt^{1}\wedge ...\wedge dt^{p}\right) ^{\frac{1}{2%
}}.
\]
The condition $\left\| \widetilde{u}\right\| \rightarrow \infty $ implies
$$
\int_{T_{0}}\left| \frac{\partial u}{\partial t}\left( t\right) \right|
^{2}dt^{1}\wedge ...\wedge dt^{p}\rightarrow \infty .\eqno(6)
$$
From the hypothesis and $\left( 5\right) $ or $\left( 6\right) $ it follows
that if $\left\| u\right\| \rightarrow \infty $, then $\varphi \left(
u\right) \rightarrow \infty $. So $\varphi $ is a coercitive application.
This means that $\varphi $ has a minimizing bounded sequence $\left(
u_{k}\right) $. The Hilbert space $H_{T}^{1}$ is reflexive. By consequence,
the sequence $\left( u_{k}\right) $ (or one of his subsequence) is weakly
convergent in $H_{T}^{1}$ with the limit $u$. Because
\[
\varphi _{2}\left( u\right) =\int_{T_{0}}\delta _{ij}\delta ^{\alpha \beta }%
\frac{\partial u^{i}}{\partial t^{\alpha }}\left( t\right) \frac{\partial
v^{j}}{\partial t^{\beta }}\left( t\right) dt^{1}\wedge ...\wedge dt^{p}
\]
is convex, it follows that $\varphi _{2}$ is weakly lower semi-continuous,
so that the action
\[
\varphi \left( u\right) =\varphi _{1}\left( u\right) +\varphi _{2}\left(
u\right)
\]
is weakly lower semi-continuous and $\varphi \left( u\right) \leq \underline{%
\lim }\varphi \left( u_{k}\right) .$ This means that $u$ is minimum point of
$\varphi .$

We build the function
\[
\Phi :[-1,1]\rightarrow R,
\]
\[
\Phi \left( \lambda \right) =\varphi \left( u+\lambda v\right)
\]
\[
=\int_{T_{0}}\left[ \frac{1}{2}\left| \frac{\partial }{\partial t}\left(
u\left( t\right) +\lambda v\left( t\right) \right) \right| ^{2}+F\left(
t,u\left( t\right) +\lambda v\left( t\right) \right) \right] dt^{1}\wedge
...\wedge dt^{p},
\]
where $v\in C_{T}^{\infty }.$ The point $\lambda =0$ is a critical point of $%
\Phi $ if and only if the point $u$ is a critical point of $\varphi $.
Consequently
\[
0=\left\langle \varphi ^{\prime }\left( u\right) ,v\right\rangle
=\int_{T_{0}}\left[ \delta ^{\alpha \beta }\delta _{ij}\frac{\partial u^{i}}{%
\partial t^{\alpha }}\frac{\partial v^{j}}{\partial t^{\beta }}+\delta
_{ij}\nabla ^{i}F\left( t,u\left( t\right) \right) v^{j}\left( t\right) %
\right] dt^{1}\wedge ...\wedge dt^{p},
\]
for all $v\in H_{T}^{1}$ and hence for all $v\in C_{T}^{\infty }.$ According
to the definition of the weak divergence, i.e.,
\[
\int_{T_{0}}\delta ^{\alpha \beta }\delta _{ij}\frac{\partial u^{i}}{%
\partial t^{\alpha }}\frac{\partial v^{j}}{\partial t^{\beta }}dt^{1}\wedge
...\wedge dt^{p} =-\int_{T_{0}}\delta ^{\alpha \beta }\delta _{ij}\frac{%
\partial ^{2}u^{i}}{\partial t^{\alpha }\partial t^{\beta }}%
v^{j}dt^{1}\wedge ...\wedge dt^{p},
\]
the Jacobi matrix function $\displaystyle\frac{\partial u}{\partial t}$ has
weak divergence (the function $u$ has a weak Laplacian) and
\[
\bigtriangleup u\left( t\right) =\nabla F\left( t,u\left( t\right) \right)
\]
a.e. on $T_{0}$. Also, the existence of weak derivatives $\displaystyle\frac{%
\partial u}{\partial t}$ and $\bigtriangleup u$ implies that
\[
u\mid _{S_{i}^{-}}=u\mid _{S_{i}^{+}},\frac{\partial u}{\partial t}\mid
_{S_{i}^{-}}=\frac{\partial u}{\partial t}\mid _{S_{i}^{+}}.
\]

{\bf Remark}. If the function $u$ is at least of class $C^2$, then the
definition of the weak divergence of the Jacobian matrix ${\frac{\partial u}{%
\partial t}}$ (or of the weak Laplacian $\bigtriangleup u$) coincides with
the classical definition. This fact is obvious if we have in mind the
formula of {\it integration by parts}

\[
\int_{T_{0}}\delta ^{\alpha \beta }\delta _{ij}\frac{\partial u^{i}}{%
\partial t^{\alpha }}\frac{\partial v^{j}}{\partial t^{\beta }}dt^{1}\wedge
...\wedge dt^{p}
\]
\[
=\int_{T_{0}}\delta ^{\alpha \beta }\delta _{ij}\frac{\partial }{\partial
t^{\alpha }}\left( \frac{\partial u^{i}}{\partial t^{\alpha }}v^{j}\right)
dt^{1}\wedge ...\wedge dt^{p}-\int_{T_{0}}\delta ^{\alpha \beta }\delta _{ij}%
\frac{\partial ^{2}u^{i}}{\partial t^{\alpha }\partial t^{\beta }}%
v^{j}dt^{1}\wedge ...\wedge dt^{p}.
\]

\bigskip

\centerline{\bf References}

\bigskip

[1] L. V. Kantorovici, G. P. Akilov: {\it Analiz\u{a} func\c{t}ional\u{a}},
Editura \c{s}tiin\c{t}ific\u{a} \c{s}i enciclopedic\u{a}, Bucure\c{s}ti,
1980.

[2] J. Mawhin, M. Willem: {\it Critical Point Theory and Hamiltonian Systems}%
, Springer-Verlag, 1989.

[3] S. G. Mihlin: {\it Ecua\c{t}ii liniare cu derivate par\c{t}iale},
Editura \c{s}tiin\c{t}ific\u{a} \c{s}i enciclopedic\u{a}, Bucure\c{s}ti,
1983.

[4] R. E. Showalter: {\it Hilbert Space Methods for Partial Differential
Equations}, Electronic Journal of Differential Equations Monograph 01,
1994.

[5] C. Udri\c{s}te, M. Postolache: {\it Atlas of Magnetic Geometric Dynamics}%
, Geometry Balkan Press, Bucharest, 2001.

[6] C. Udri\c{s}te: {\it From Integral Manifolds and Metrics to Potential
Maps}, Conference Michigan State University, April 13-27, 2001; Eleventh
Midwest Geometry Conference, Wichita State University, April 27-29, 2001,
Atti del' Academia Peloritana dei Pericolanti, Clase 1 di Scienze Fis. Mat.
e Nat., 81-82, A 01006 (2003-2004), 1-14.

[7] C. Udri\c{s}te, I. Duca: {\it Periodical Solutions of Multi-Time
Hamilton Equations,} Analele Universita\c{t}ii Bucure\c{s}ti, 55, 1 (2005),
178-188.

[8] C. Udri\c{s}te, M. Ferrara, D. Opri\c{s}: {\it Economic Geometric
Dynamics, }Geometry Balkan Press, Bucharest, 2004.

[9] C. Udri\c{s}te, A-M. Teleman: {\it Hamilton Approaches of Fields Theory}%
, IJMMS, 57 (2004), 3045-3056; ICM Satelite Conference in Algebra and
Related Topics, University of Hong-Kong, 13-18.08.02.

[10] C. Udri\c{s}te, M. Neagu:{\it \ From PDE Systems and Metrics to
Generalized Field Theories}, http://xxx.lanl.gov/abs/math.DG/0101207.

[11] C. Udri\c{s}te: {\it Nonclassical Lagrangian Dynamics and Potential Maps%
}, The Conference in Mathematics in Honour of Professor Radu Ro\c{s}ca on
the Occasion of his Ninetieth Birthday, Katholieke University Brussel,
Katholieke University Leuven, Belgium, Dec.11-16, 1999;

http://xxx.lanl.gov/math.DS/0007060.

[12] C. Udri\c{s}te: {\it \ Solutions of DEs and PDEs as Potential Maps
Using First Order Lagrangians}, Centenial Vr\^{a}nceanu, Romanian Academy,
University of Bucharest, June 30-July 4, (2000); Balkan Journal of Geometry and Its
Applications 6, 1, 93-108, 2001; http://xxx.lanl.gov/math.DS/0007061.

\end{document}